\documentclass[oneside,11pt,reqno]{amsart}

\usepackage[T1]{fontenc}
\usepackage{lmodern}
\usepackage{microtype}
\usepackage{xcolor}
\usepackage[margin=1in]{geometry}
\usepackage{amsmath,amssymb,amsthm,mathtools,booktabs}

\usepackage{enumitem}

\newtheorem{theorem}{Theorem}[section]
\newtheorem{corollary}[theorem]{Corollary}
\newtheorem{lemma}[theorem]{Lemma}
\newtheorem{proposition}[theorem]{Proposition}
\theoremstyle{definition}

\newtheorem{example}[theorem]{Example}

\theoremstyle{remark}

\DeclareMathOperator{\br}{br}
\numberwithin{equation}{section}
\numberwithin{figure}{section}
\numberwithin{table}{section}

\title{Equal knapsack identities between symmetric group character degrees}

\author{David J.\ Hemmer}
\address{Department of Mathematical Sciences\\
Michigan Technological University\\
Houghton, MI 49931}
\email{djhemmer@mtu.edu}

\author{Armin Straub}
\address{Department of Mathematics and Statistics\\
University of South Alabama\\
Mobile, AL 36688}
\email{straub@southalabama.edu}

\author{Karlee J.\ Westrem}
\address{Department of Mathematical Sciences\\
Appalachian State University\\
Boone, NC 28608}
\email{westremk@appstate.edu}

\usepackage[hidelinks]{hyperref}
\dedicatory{Dedicated to the memory of Jonathan Alperin.}

\thanks{We are honored to dedicate this paper to the memory of Jonathan Alperin, who advised the first author and was the mathematical grandfather of the third. Jon always loved the representation theory of the symmetric group. He used to call it his ``laboratory''. Indeed he and Fong proved the famous Alperin weight conjecture in the case of symmetric groups shortly after the general conjecture was announced. He also liked to say that if you proved something about symmetric groups that it was guaranteed to last for centuries, in contrast to a result about other more obscure algebraic objects that might rise and fall in popularity over time. The ``symmetric group is forever'' he would say. Jon was a great mathematician and a great advisor, and he will be sorely missed.}

\begin{document}

\begin{abstract}
We prove a series of ``knapsack'' type equalities for irreducible character degrees of symmetric groups. That is, we find disjoint subsets of the partitions of $n$ so that the two corresponding character-degree sums are equal. Our main result refines our recent description of the Riordan numbers as the sum of all character degrees $f^\lambda$ where $\lambda$ is a partition of $n$ into three parts of the same parity. In particular, the sum of the ``fat-hook'' degrees $f^{(k,k,1^{n-2k})}+f^{(k+1,k+1,1^{n-2k-2})}$ equals the sum of all $f^\lambda$ where $\lambda$ has three parts, with the second equal to $k$ and the second and third of equal parity. We further prove an infinite family of additional ``knapsack'' identities between character degrees.
\end{abstract}

\maketitle
\section{Introduction}
\label{sec: Introduction}
Let $\Sigma_n$ be the symmetric group on $n$ letters and recall that the complex irreducible characters of $\Sigma_n$ are labeled by partitions of $n$. The degree of the irreducible character corresponding to a partition $\lambda \vdash n$ is the number of standard Young tableaux of shape $\lambda$, which we denote by $f^\lambda$.

The famous Catalan numbers $C(n)$ enumerate hundreds of different sets \cite{StanleyCatalan}, including Dyck paths of semilength $n$ (that is, paths from $(0,0)$ to $(2n,0)$ using steps $U = (1,1)$ and $D = (1,-1)$ and not going below the $x$-axis).
They also count standard Young tableaux of shape $(n,n)$, that is $C(n)=f^{(n,n)}$. Allowing also horizontal steps, one obtains \emph{Motzkin paths}, which are paths from $(0,0)$ to $(n,0)$ using only steps $U = (1,1)$, $F = (1,0)$, and $D = (1,-1)$, and not going below the $x$-axis. These are enumerated by the Motzkin numbers $M(n)$, sequence A001006 in \cite{oeis}. $M(n)$ is also well known to count the number of standard Young tableaux of size $n$ with three or fewer rows; for an example of a bijection, see \cite{MatsakisMotzkinInspired}. In our notation, this count can be expressed as
\begin{equation}
\label{eq:Motzkinnumberthreepartdegrees}
M(n)=\sum_{\lambda=(\lambda_1,\lambda_2,\lambda_3) \vdash n}f^\lambda.
\end{equation}
Here, and elsewhere, we allow parts to be zero. As such, \eqref{eq:Motzkinnumberthreepartdegrees} sums over all partitions with up to three nonzero parts.

A \emph{Riordan path} is a Motzkin path with the additional requirement that there may not be a flat step $F$ on the $x$-axis. We let $R(n)$ be the number of Riordan paths of length $n$; this is the sequence A005043 in \cite{oeis}. The Riordan numbers $R(n)$ also have the following interpretation in terms of degrees of irreducible symmetric group characters. This is a comment of Regev given on the OEIS entry (for a proof see \cite{HemmerStraubWestrem2025}).

\begin{proposition}\cite[Proposition 1.2]{HemmerStraubWestrem2025}
  \label{prop:RiordannumberSYT}
  Let $0 \leq m < n$. The number of Riordan paths of length $n$ with $m$ flat steps and $k$ up steps (and thus $k$ down steps) is $f^{(k,k,1^m)}$.
\end{proposition}

In \cite{HemmerStraubWestrem2025}, we gave an interpretation of $R(n)$ as the sum of the $f^\lambda$, not over all three-part partitions like $M(n)$, but instead over those three-part partitions where the three parts have equal parity (necessarily the parity of $n$).

\begin{theorem}\cite[Theorem 4.5]{HemmerStraubWestrem2025}
  \label{thm:equalitySpecht}
  Let $X = \{ (\lambda_1, \lambda_2, \lambda_3) \vdash n \mid \lambda_1 \equiv \lambda_2 \equiv \lambda_3 \pmod{2} \}$ and let $Y = \{ (k,k,1^{n-2k}) \mid 1 \leq k \leq \lfloor n/2 \rfloor \}$. Then:
  \begin{equation}
    \sum_{\lambda \in X} f^\lambda = \sum_{\mu \in Y} f^\mu=R(n).
    \label{eq: SumsSpechtequal}
  \end{equation}
\end{theorem}

We note that the bijection in \cite{MatsakisMotzkinInspired} proving \eqref{eq:Motzkinnumberthreepartdegrees} does not send Riordan paths to tableaux with shapes with parts of equal parity. Indeed, we do not know \footnote{\label{fn} Note added in proof: We thank the referee for pointing out that such a bijection may be found in  \cite{Jagenteufel}.} of a bijective proof of Theorem \ref{thm:equalitySpecht}.

In this paper, we prove that the equality in \eqref{eq: SumsSpechtequal} may be refined to a series of equalities, one for each choice of the second part $\lambda_2$. For example, when $n=20$ we prove in \eqref{eq:firstdimension}:
\begin{align*}
f^{(20)} &= f^{(1^{20})} \\
f^{(18,2)} + f^{(16,2,2)} &= f^{(2,2,1^{16})} + f^{(3,3,1^{14})} \\
f^{(16,4)} + f^{(14,4,2)} + f^{(12,4,4)} &= f^{(4,4,1^{12})} + f^{(5,5,1^{10})} \\
f^{(14,6)} + f^{(12,6,2)} + f^{(10,6,4)} + f^{(8,6,6)}
  &= f^{(6,6,1^{8})} + f^{(7,7,1^{6})} \\
f^{(12,8)} + f^{(10,8,2)} + f^{(8,8,4)} &= f^{(8,8,1^{4})} + f^{(9,9,1^{2})} \\
f^{(10,10)} &= f^{(10,10)}.
\end{align*}

We attempted to prove these identities by induction using the branching rule for symmetric group irreducible characters. However, the branching rule does not preserve the parity conditions, which required us to discover and prove three additional sets of identities all together using a single induction proof. Yet, even these four identities together are not preserved under the branching rules, so attempting a proof by induction leads to error terms which are not subject to the inductive hypothesis. These error terms led to the discovery of two more sets of identities (see Section~\ref{sec:auxiliary}) between character degrees which may be proved directly using the hook length formula, and thus let us prove the four main identities by induction.

These four identities are recorded in the following theorem which is our first main result. Note that the identity \eqref{eq:firstdimension}, which always holds when $k \equiv n \pmod{2}$, is sufficient for the refinement of Theorem~\ref{thm:equalitySpecht} which we described above.

\begin{theorem}\label{thm: Generalizedknapsack}
  Let $n$ and $k$ be positive integers such that $k \le \tfrac{n}{2} - 1$.
  Let
  $$X_0(n,k) = \{ (\lambda_1, k, \lambda_3)\vdash n \mid  k \equiv  \lambda_3 \pmod{2}\}$$
  and
  $$X_1(n,k) = \{ (\lambda_1, k, \lambda_3)\vdash n \mid  k \not\equiv  \lambda_3 \pmod{2} \}.$$
  If $k \leq \lceil\frac{n}{3}\rceil$ or if $k \equiv n \pmod{2}$ then
  \begin{equation}
    \label{eq:firstdimension}
    \sum_{\lambda \in X_0(n,k)} f^\lambda = f^{(k,k, 1^{n-2k})} + f^{(k+1,k+1,1^{n-2k-2})}
  \end{equation}
  as well as
  \begin{equation}
    \label{eq:seconddimension}
    \sum_{\lambda \in X_1(n,k)} f^\lambda = f^{(k+1,k,1^{n-2k-1})}.
  \end{equation}
  If $k > \lceil\frac{n}{3}\rceil$ and $k \not \equiv n \pmod{2}$ then the equalities above hold with the roles of $X_0$ and $X_1$ switched.
\end{theorem}

If we adopt the convention that any term $f^\lambda$ on the right-hand sides of \eqref{eq:firstdimension} and \eqref{eq:seconddimension} is zero if $\lambda$ is not a partition, then the statement of Theorem~\ref{thm: Generalizedknapsack} continues to hold for all $n\geq 3$ and $k\geq 0$.

Above, we already illustrated equation \eqref{eq:firstdimension} when $n=20$ and $k \equiv n \pmod{2}$. In the following, we briefly illustrate both of the equalities \eqref{eq:firstdimension} and \eqref{eq:seconddimension} for some additional choices of $n$ and $k$.  For instance, for $n=32$, $k=11$, the two equalities become
\begin{align*}
f^{(20,11,1)}+f^{(18,11,3)}+f^{(16,11,5)}+f^{(14,11,7)}+f^{(12,11,9)}
&= f^{(11,11,1^{10})}+f^{(12,12,1^8)},\\
f^{(21,11,0)}+f^{(19,11,2)}+f^{(17,11,4)}+f^{(15,11,6)}+f^{(13,11,8)}+f^{(11,11,10)}
&= f^{(12,11,1^9)}.
\end{align*}
Similarly, for $n=32$, $k=12$, the equalities are
\begin{align*}
f^{(20,12,0)}+f^{(18,12,2)}+f^{(16,12,4)}+f^{(14,12,6)}+f^{(12,12,8)}
&=f^{(12,12,1^{8})}+f^{(13,13,1^6)},\\
f^{(19,12,1)}+f^{(17,12,3)}+f^{(15,12,5)}+f^{(13,12,7)}
&=f^{(13,12,1^{7})}.
\end{align*}
On the other hand, for $n=32$, $k=13$, we obtain
\begin{align*}
f^{(18,13,1)}+f^{(16,13,3)}+f^{(14,13,5)}
&=f^{(14,13,1^5)},\\
f^{(19,13,0)}+f^{(17,13,2)}+f^{(15,13,4)}+f^{(13,13,6)}&=f^{(13,13,1^6)}+f^{(14,14,1^4)}.
\end{align*}

In the next section, we prove the two auxiliary identities required in the inductive proof of Theorem~\ref{thm: Generalizedknapsack}, and then in Section \ref{sec: Proof of main theorem} we outline the proof of  Theorem \ref{thm: Generalizedknapsack}, with full details given in the appendix. In Section~\ref{sec: Extension}, we show that the auxiliary identities of Section~\ref{sec:auxiliary} can be considerably extended.  As an application of these extensions, which are our second main result, we give an alternative proof of identity \eqref{eq:firstdimension} of Theorem~\ref{thm: Generalizedknapsack}. We conclude by posing additional problems and conjecturing further identities.

\section{Auxiliary Identities}
\label{sec:auxiliary}

In this section we prove two identities between character degrees which appear of independent interest and are a necessary ingredient to proving Theorem~\ref{thm: Generalizedknapsack} by induction using the branching rule. Unlike the identities in Theorem~\ref{thm: Generalizedknapsack}, the two auxiliary identities are of fixed length and so we are able to prove them directly from the hook length formula.

Recall that, for a partition $\lambda\vdash n$ with conjugate partition $\lambda'$, the number of standard Young tableaux of shape $\lambda$ is given by the famous \emph{hook length formula} as
\begin{equation}
\label{eq:hook-length-formula}
 f^\lambda = \frac{n!}{\prod h_{i,j}}.
\end{equation}
where the product is over the boxes in the Young diagram for $\lambda$ and $h_{i,j}=\lambda_i+\lambda_j'-i-j+1$ is the hook length at position $(i,j)$.

Using the hook length formula one readily obtains the following explicit formulas.
The proofs are standard and rely on the fact that for diagrams of shape $(a,b,1^t)$ or $(r,s,t)$ it is easy to write down all the hook lengths.
For completeness, we include details for the second formula. The proof of the first formula is similar and, hence, omitted.

\begin{lemma}
\label{lem:closedformab111}
For $a\ge b\ge1$, $t\ge0$,
\[
  f^{(a,b,1^t)}
  \;=\;
  \frac{(a+b+t)!\,(a-b+1)}{(a+t+1)(b+t)\;a!\,(b-1)!\;t!}.
\]
\end{lemma}

\begin{lemma}\label{lem:closedformrst}
For $r\ge s\ge t\ge0$,
\[
  f^{(r,s,t)}
  \;=\;
  \frac{(r+s+t)!(r-s+1)(r-t+2)(s-t+1)}{(r+2)!(s+1)!\,t!}.
\]
\end{lemma}

\begin{proof}
We simply need to determine the hook lengths contributing to the denominator
in \eqref{eq:hook-length-formula}. For the third row these are $\{1,2,\dots,
t\}$ which contributes $t!$ to the denominator. For the second row they are
$$\{1,2,\dots,s-t,s-t+2,s-t+3,\dots,s+1\}$$ which contributes
$\frac{(s+1)!}{s-t+1}$.  And finally for the top row the hook lengths are
$$\{1,2,\dots,r-s,r-s+2,r-s+3,\dots,r-t+1,r-t+3,r-t+4,\dots,r+2\}$$ which
contributes $\frac{(r+2)!}{(r-s+1)(r-t+2)}$.
Combined, the claimed formula thus follows from \eqref{eq:hook-length-formula}.
\end{proof}

The identities in Theorem~\ref{thm: Generalizedknapsack} are relationships between various $f^\lambda$ for partitions $\lambda$ of these shapes. The following auxiliary identity is a fixed-length relationship of that kind as well.
Fix integers $k\ge1$ and $m\ge4$ and define
\begin{equation}
\label{eq: DefineL(k,m)}
  L(k,m)= f^{(k,k,1^m)}+f^{(k+1,k+1,1^{m-2})}+f^{(k+2,k+2,1^{m-4})}.
\end{equation}
Then we have:

\begin{lemma}
\label{lem:3=2identity}
For $k\ge1$, $m\ge4$,
$$
  L(k,m)
  = f^{(k+2,k,1^{m-2})}
  +
  \begin{cases}
  f^{(k,k,m)}, & m\le k,\\
  0, & m=k+1\ \text{or}\ m=k+2,\\
  f^{(m-2,k+1,k+1)}, & m\ge k+3.
  \end{cases}
$$
\end{lemma}

The following three examples, all for $n=35$, illustrate each of the three cases from Lemma~\ref{lem:3=2identity}:
\begin{align*}
f^{(14,14,1^7)} + f^{(15,15,1^5)} + f^{(16,16,1^3)} &= f^{(16,14,1^5)} + f^{(14,14,7)}, \\
f^{(11,11,1^{13})} + f^{(12,12,1^{11})} + f^{(13,13,1^9)} &= f^{(13,11,1^{11})},\\
f^{(7,7,1^{21})} + f^{(8,8,1^{19})} + f^{(9,9,1^{17})} &= f^{(9,7,1^{19})} + f^{(19,8,8)}
.
\end{align*}

\begin{proof}
  Consider the difference
  \begin{equation*}
    \Delta_3 (k, m) = f^{(k, k, 1^m)} + f^{(k + 1, k + 1, 1^{m - 2})} +
     f^{(k + 2, k + 2, 1^{m - 4})} - f^{(k + 2, k, 1^{m - 2})}
  \end{equation*}
  and apply Lemma~\ref{lem:closedformab111} to each term. We obtain, in
  particular,
  \begin{eqnarray}
    f^{(k, k, 1^m)} & = & \frac{(2 k + m) !}{k! (k - 1) !m!} \, \frac{1}{(k +
    m) (k + m + 1)},  \label{eq:f:kk1m}\\
    f^{(k + 2, k, 1^{m - 2})} & = & \frac{(2 k + m) !}{(k + 2) ! (k - 1) ! (m
    - 2) !} \, \frac{3}{(k + m + 1) (k + m - 2)} . \nonumber
  \end{eqnarray}
  For each of the four terms $f^{\lambda}$ of $\Delta_3 (k, m)$, we find that
  $f^{\lambda} / f^{(k, k, 1^m)}$ is a rational function in $k$ and $m$. For
  instance,
  \begin{equation*}
    \frac{f^{(k + 2, k, 1^{m - 2})}}{f^{(k, k, 1^m)}} = \frac{3 (k + m) (m -
     1) m}{(k + 1) (k + 2) (k + m - 2)} .
  \end{equation*}
  Consequently, $\Delta_3 (k, m) / f^{(k, k, 1^m)}$ is a rational function in
  $k$ and $m$ as well. Indeed, this rational function factors as
  \begin{equation*}
    \frac{\Delta_3 (k, m)}{f^{(k, k, 1^m)}} = \frac{(k - m + 1) (k - m + 2)
     (k + m) (k + m + 1)}{k (k + 1)^2 (k + 2)} .
  \end{equation*}
  Multiplying both sides with $f^{(k, k, 1^m)}$ as evaluated in
  \eqref{eq:f:kk1m}, we obtain
  \begin{equation}
    \Delta_3 (k, m) = \frac{(2 k + m) !}{(k + 1) ! (k + 2) !m!} (k - m + 1) (k
    - m + 2) . \label{eq:D3:hyp}
  \end{equation}
  This shows that $\Delta_3 (k, m) = 0$ if $k = m - 1$ or $k = m - 2$, as
  claimed. 
  
  On the other hand, by the hook length formula as in
Lemma~\ref{lem:closedformrst} we have
\[
  f^{(m - 2, k + 1, k + 1)}
  =
  \frac{(2 k + m) ! (m - k - 2) (m - k - 1)}
       {m! (k + 2)! (k + 1)!}
\]
provided that $m - 2 \geq k + 1$, that is, provided that
$m \geq k+3$. Since the right-hand side equals
\eqref{eq:D3:hyp}, this proves the case $m \geq k+3$.
Likewise applying Lemma~\ref{lem:closedformrst}, we find that, when
$m\leq k$,
\[
  f^{(k,k,m)}
  =
  \frac{(2k+m)!(k-m+1)(k-m+2)}
       {(k+2)!(k+1)!\,m!}.
\]
This is precisely the expression in \eqref{eq:D3:hyp}. This proves the case $m\leq k$. Finally, the
remaining cases are exactly $m=k+1$ and $m=k+2$, where
\eqref{eq:D3:hyp} gives $\Delta_3(k,m)=0$.
\end{proof}

Theorem \ref{thm: Generalizedknapsack} splits into cases depending on how the size of $k$ compares to $n/3$ and the respective parities of $k$ and $n$.  When $k$ is almost equal to $n/3$ we get a sort of boundary case in the induction, which needs the following identity.

\begin{proposition}
\label{prop: boundaryspecialcase}
Let $k$ and $m$ be positive integers with $m \ge 2$. If $k= m \pm 1$, then
$$f^{(k,k,1^m)}+f^{(k{+}1,k{+}1,1^{m-2})}=f^{(k{+}1,k,1^{m-1})}.$$
\end{proposition}

\begin{proof}
One readily verifies from Lemma~\ref{lem:closedformab111} that
\begin{equation}
\label{eq: special boundary case}
f^{(k,k,1^m)}+f^{(k{+}1,k{+}1,1^{m-2})}-f^{(k{+}1,k,1^{m-1})}
=\frac{(2k+m)!\,(k-m-1)(k-m+1)}{(k+m-1)(k+m+1)\,(k+1)!\,k!\,m!}.
\end{equation}
The claimed identity then follows because the right-hand side of \eqref{eq: special boundary case} is zero when $k=m\pm 1$.
\end{proof}

\section{Proof of Theorem \ref{thm: Generalizedknapsack} via the branching rule}
\label{sec: Proof of main theorem}
Recall that the branching rule describes how irreducible characters of the symmetric group $\Sigma_n$ decompose upon restriction to $\Sigma_{n-1}$. For $\lambda \vdash n$, the dimension $f^\lambda$ is equal to the sum of all the $f^\mu$ where $\mu \vdash n-1$ is obtained by removing a box from the Young diagram of $\lambda$ in a way that leaves a valid partition of $n-1$.

The proof of Theorem \ref{thm: Generalizedknapsack} is by induction on $n$. Since small values of $n$ are easily checked, we assume that all possible equalities from the theorem hold for all values smaller than $n$. We write down each of the desired equalities and apply the branching rule on both sides. On the right-hand side we will get four terms from \eqref{eq:firstdimension} or three terms from \eqref{eq:seconddimension}. On the left-hand side we look at the partitions in $X_0(n,k)$ or $X_1(n,k)$, apply the branching theorem, and divide the result into three summands based on whether boxes were removed from the first, second or third rows.

About half the time all three rows give sets of the form $X_i(n-1,k)$ (rows 1 and 3) or $X_i(n-1,k-1)$ (row 2) for $i=0$ or $1$, so the inductive hypotheses apply to all three terms and give the result immediately. The other half of the time there is a term of the form $f^{(a,a,b)}$ or $f^{(a,b,b)}$ that is missing from one of the resulting sums, and we need to apply Lemma~\ref{lem:3=2identity}. In one boundary case when $3 \mid n$ and $k=\frac{n}{3}+1$ we also need Proposition \ref{prop: boundaryspecialcase}.

Let's begin by looking at some illustrative examples. Suppose $n=20$ and $k=5$. To prove \eqref{eq:firstdimension} we need to show that
\begin{equation}
\label{eq: firstexamplemaintheorem}
f^{(14,5,1)}+f^{(12,5,3)}+f^{(10,5,5)}=f^{(5,5,1^{10})}+f^{(6,6,1^8)}.
\end{equation}
Applying the branching rule on the right-hand side we get
$$f^{(5,4,1^{10})}+f^{(5,5,1^9)}+f^{(6,5,1^8)}+f^{(6,6,1^7)}.$$
To likewise apply the branching rule on the left-hand side of \eqref{eq: firstexamplemaintheorem}, we remove boxes from one row at a time from the involved partitions. Removing from the first row, we obtain $f^{(13,5,1)}+f^{(11,5,3)}+f^{(9,5,5)}$ which equals $f^{(5,5,1^9)}+f^{(6,6,1^7)}$ by the induction hypothesis. Removing instead from the second row, we get $f^{(14,4,1)}+f^{(12,4,3)}$ which equals $f^{(5,4,1^{10})}$ by induction. Finally, removing from the third row yields $f^{(14,5)}+f^{(12,5,2)}+f^{(10,5,4)}$ which is $f^{(6,5,1^8)}$ by induction. The total from the left-hand side also is $f^{(5,4,1^{10})}+f^{(5,5,1^9)}+f^{(6,5,1^8)}+f^{(6,6,1^7)}$, thus proving the identity \eqref{eq: firstexamplemaintheorem}.

To likewise prove \eqref{eq:seconddimension} for the same $n$ and $k$, we need to show that
\begin{equation}
\label{eq: secondexamplemaintheorem}
f^{(15,5)}+f^{(13,5,2)}+f^{(11,5,4)}=f^{(6,5,1^9)}.
\end{equation}
Applying the branching rule on the right-hand side we get
$$f^{(5,5,1^9)}+f^{(6,4,1^9)}+f^{(6,5,1^8)}.$$
To apply the branching rule on the left-hand side, we again remove boxes from one row at a time. Removing from the first row results in $f^{(14,5)}+f^{(12,5,2)}+f^{(10,5,4)}$ which equals $f^{(6,5,1^8)}$ by the induction hypothesis. Removing instead from the second row, we obtain $f^{(15,4)}+f^{(13,4,2)}+f^{(11,4,4)}$ which equals $f^{(4,4,1^{11})}+f^{(5,5,1^{9})}$ by induction. Removing from the third row, we get $f^{(13,5,1)}+f^{(11,5,3)}$. The inductive hypothesis for $X_0(19,5)$ gives a formula for $f^{(13,5,1)}+f^{(11,5,3)}+f^{(9,5,5)}$. The partition $(9,5,5)$ lies in $X_0(19,5)$ but does not show up in the branching (it would have to come from $(9,5,6)$ which is not a partition.)  By the induction hypothesis, the terms resulting from removing from the third row therefore are
$$f^{(5,5,1^9)}+f^{(6,6,1^7)}-f^{(9,5,5)}.$$
We are left to confirm that
$$f^{(5,5,1^9)}+f^{(6,4,1^9)}+f^{(6,5,1^8)}=f^{(6,5,1^8)}+f^{(4,4,1^{11})}+f^{(5,5,1^{9})}+f^{(5,5,1^9)}+f^{(6,6,1^7)}-f^{(9,5,5)}.$$
Cancelling the terms that appear on both sides and applying Lemma~\ref{lem:3=2identity} for $k=4, m=11$, we readily confirm this equality, thus completing the proof of \eqref{eq:seconddimension} in the present example.

Suppose we want to instead prove \eqref{eq:seconddimension} for $n=14$ and $k=5$ (so that $k = \lceil n/3 \rceil$). This is a ``boundary case'' and we find that the expected ``missing term'' does not appear. In this case, we  need to show that
$$f^{(9,5)}+f^{(7,5,2)}+f^{(5,5,4)}=f^{(6,5,1^3)}.$$
When we remove from the third rows on the left-hand side, we get $f^{(7,5,1)}+f^{(5,5,3)}$. This matches $X_0(13,5)=\{(7,5,1),(5,5,3)\}$ so that there is not a missing term of the form $f^{(a,5,5)}$ (for contrast, note that for $n=16$ the set $X_0(15,5)$ contains $(5,5,5)$ and this partition is missing in the branching). Continuing to branch as in the previous examples, we end up needing to prove that
$$f^{(5,5,1^3)}+f^{(6,4,1^3)}+f^{(6,5,1^2)}=f^{(6,5,1^2)}+f^{(4,4,1^5)}+f^{(5,5,1^3)}+f^{(5,5,1^3)}+f^{(6,6,1)}.$$
This reduces to
$$f^{(6,4,1^3)}=f^{(4,4,1^5)}+f^{(5,5,1^3)}+f^{(6,6,1)},$$
which is exactly the case $k=4, m=5$ where the second term in Lemma~\ref{lem:3=2identity} vanishes. This happens precisely when that term is not needed to cancel anything out.  This illustrates how the three cases in Lemma~\ref{lem:3=2identity} are precisely calibrated to the needs of our inductive proof.

Finally, there is a boundary case where $n=3x$ and $k=x+1$ (so that $k > \lceil\frac{n}{3}\rceil=x$). This is the ``large $k$'' regime (meaning that $k > \lceil\frac{n}{3}\rceil$) of Theorem~\ref{thm: Generalizedknapsack} where $k \not \equiv n \pmod{2}$. On the other hand, when we branch by removing from the second row, the value of $k$ drops to $x$ and $\lceil\frac{n-1}{3}\rceil=x$ which is a switch to the ``small $k$'' regime (meaning that $k \leq \lceil\frac{n}{3}\rceil$).  In this case, we need Proposition~\ref{prop: boundaryspecialcase}. For example, suppose $n=24$ and $k=9$ and we try to prove \eqref{eq:firstdimension}, which is
$$f^{(14,9,1)}+f^{(12,9,3)}+f^{(10,9,5)}=f^{(10,9,1^5)}.$$
Applying the branching rule to both sides as before and cancelling equal terms, we end up having to prove that
$$f^{(10,10,1^3)}+f^{(9,8,1^6)}=f^{(10,8,1^5)}+f^{(8,8,7)}.$$
By Lemma~\ref{lem:3=2identity} the right-hand side is just $f^{(8,8,1^7)}+f^{(9,9,1^5)}+f^{(10,10,1^3)}$ so we end up needing to show that
$$f^{(9,8,1^6)}=f^{(8,8,1^7)}+f^{(9,9,1^5)},$$
which is an instance of Proposition \ref{prop: boundaryspecialcase}.

At this point, we could safely leave the general details to the reader. Theorem \ref{thm: Generalizedknapsack} contains two equalities which are flipped when $k > \lceil\frac{n}{3}\rceil$ and $k \not \equiv n \pmod{2}$. So one has to consider the cases where $k$ is quite large or quite small compared to $n/3$, where it is the same or opposite parity as $n$, and for both equalities. Then, as we indicated above, there are boundary cases when $k$ is close to $n/3$ where the middle case of Lemma~\ref{lem:3=2identity} or Proposition \ref{prop: boundaryspecialcase} are used.

The details of all these cases may be found in Appendix~\ref{sec:complete-branching-proof}. On the other hand, we provide in the next section an alternative analytical proof of the part of Theorem~\ref{thm: Generalizedknapsack} that refines Theorem~\ref{thm:equalitySpecht}.

\section{An extension and an alternative proof}\label{sec: Extension}

Let us consider the sums
\begin{equation*}
  L_d (k, m) = \sum_{j = 0}^{d - 1} f^{(k + j, k + j, 1^{m - 2 j})} .
\end{equation*}
In terms of this notation, Lemma~\ref{lem:3=2identity} shows that
\begin{equation*}
  L_3 (k, m) = f^{(k + 2, k, 1^{m - 2})} + \left\{\begin{array}{ll}
     f^{(k, k, m)}, & \text{if $m \leq k$},\\
     0, & \text{if $m = k + 1$ or $m = k + 2$,}\\
     f^{(m - 2, k + 1, k + 1)}, & \text{if $m \geq k + 3$} .
   \end{array}\right.
\end{equation*}
In this section, we show that this can be extended to the sums $L_d (k, m)$
for all odd $d$. We then demonstrate that this extension provides an
alternative approach to the refinement of Theorem~\ref{thm:equalitySpecht}
proved in Theorem~\ref{thm: Generalizedknapsack}.

Before turning to the general case, we highlight the case $d = 5$.

\begin{example}
  \label{eg:L5}Let $k \geq 2$ and $m \geq 8$. If $m \leq k$
  then
  \begin{equation}
    L_5 (k, m) = f^{(k + 4, k, 1^{m - 4})} + f^{(k, k, m)} + f^{(k + 2, k, m -
    2)} + f^{(k + 2, k + 2, m - 4)} . \label{eq:L5:1}
  \end{equation}
  On the other hand, if $m \geq k + 9$ then
  \begin{equation}
    L_5 (k, m) = f^{(k + 4, k, 1^{m - 4})} + f^{(m - 2, k + 1, k + 1)} + f^{(m
    - 4, k + 3, k + 1)} + f^{(m - 6, k + 3, k + 3)} . \label{eq:L5:2}
  \end{equation}
  For each of the intermediate cases $m = k + \delta$ with $\delta \in \{ 1,
  2, \ldots, 8 \}$ a suitable variation of the identities \eqref{eq:L5:1} and
  \eqref{eq:L5:2} holds. Specifically:
  \begin{eqnarray*}
    L_5 (k, k + 1) & = & f^{(k + 4, k, 1^{k - 3})} + f^{(k + 2, k, k - 1)} +
    f^{(k + 2, k + 2, k - 3)}\\
    L_5 (k, k + 2) & = & f^{(k + 4, k, 1^{k - 2})} + f^{(k + 2, k, k)} + f^{(k
    + 2, k + 2, k - 2)}\\
    L_5 (k, k + 3) & = & f^{(k + 4, k, 1^{k - 1})} + f^{(k + 1, k + 1, k + 1)}
    + f^{(k + 2, k + 2, k - 1)}\\
    L_5 (k, k + 4) & = & f^{(k + 4, k, 1^k)} + f^{(k + 2, k + 2, k)}\\
    L_5 (k, k + 5) & = & f^{(k + 4, k, 1^{k + 1})} + f^{(k + 3, k + 1, k +
    1)}\\
    L_5 (k, k + 6) & = & f^{(k + 4, k, 1^{k + 2})} + f^{(k + 4, k + 1, k + 1)}
    + f^{(k + 2, k + 2, k + 2)}\\
    L_5 (k, k + 7) & = & f^{(k + 4, k, 1^{k + 3})} + f^{(k + 5, k + 1, k + 1)}
    + f^{(k + 3, k + 3, k + 1)}\\
    L_5 (k, k + 8) & = & f^{(k + 4, k, 1^{k + 4})} + f^{(k + 6, k + 1, k + 1)}
    + f^{(k + 4, k + 3, k + 1)}
  \end{eqnarray*}
  Note that the cases $\delta = 3, 6$ are combining terms from both of the
  identities \eqref{eq:L5:1} and \eqref{eq:L5:2} in the theorem (whereas the
  other cases result from one of the identities by removing those terms
  $f^{\lambda}$ for which $\lambda$ is not a partition).
\end{example}

The following generalizes the case $m \leq k$ of
Lemma~\ref{lem:3=2identity}. The corresponding case of large $m$ is obtained
in Corollary~\ref{cor:L:odd:sum:2}.

\begin{theorem}
  \label{thm:L:odd:sum}Let $d \geq 0$ be an integer. Let $k, m \geq
  2$ be integers such that $m \geq 4 d$. Then, for $m \leq k$,
  \begin{equation}
    L_{2 d + 1} (k, m) = f^{(k + 2 d, k, 1^{m - 2 d})} + \sum_{r = 0}^{d - 1}
    \sum_{j = 0}^r f^{(k + 2 r, k + 2 j, m - 2 (r + j))} .
    \label{eq:L:odd:sum}
  \end{equation}
\end{theorem}

\begin{proof}
  Define the functions
  \begin{eqnarray*}
    h^{(x, y, z)} & = & \frac{(x + y + z) !}{(x + 2) ! (y + 1) !z!} \; (x - y
    + 1) (x - z + 2) (y - z + 1),\\
    h_1^{(x, y ; z)} & = & \frac{(x + y + z) !}{x! (y - 1) !z!} \; \frac{(x -
    y + 1)}{(x + z + 1) (y + z)} .
  \end{eqnarray*}
  Apart from singularities, these are defined for any complex value of $x, y,
  z$ (here, as usual, the factorial $x!$ is defined as $\Gamma (x + 1)$ where
  $\Gamma$ is the gamma function). On the other hand, note that, up to
  labelling of the variables, these match the expressions obtained in
  Lemmas~\ref{lem:closedformab111} and \ref{lem:closedformrst} using the hook
  length formula. In other words, if $\lambda = (\lambda_1, \lambda_2,
  \lambda_3)$ is a partition then we have
  \begin{equation*}
    f^{(\lambda_1, \lambda_2, \lambda_3)} = h^{(\lambda_1, \lambda_2,
     \lambda_3)}
  \end{equation*}
  and similarly
  \begin{equation*}
    f^{(\lambda_1, \lambda_2, 1^r)} = h_1^{(\lambda_1, \lambda_2 ; r)} .
  \end{equation*}
  We will prove the identity~\eqref{eq:L:odd:sum} by establishing the more
  general analytic identity
  \begin{equation}
    \sum_{j = 0}^{2 d} h_1^{(k + j, k + j ; m - 2 j)} = h_1^{(k + 2 d, k ; m -
    2 d)} + \sum_{r = 0}^{d - 1} \sum_{j = 0}^r h^{(k + 2 r, k + 2 j, m - 2 (r
    + j))} \label{eq:L:odd:sum:hyp}
  \end{equation}
  where $d \geq 0$ is an integer and where $k$ and $m$ are complex
  variables. Note that the conditions on $k$ and $m$ in
  Theorem~\ref{thm:L:odd:sum} are such that each of the terms involving $h$ or
  $h_1$ in \eqref{eq:L:odd:sum:hyp} can be interpreted as $f^{\lambda}$ for
  some partition $\lambda$.
  
  The case $d = 0$ of \eqref{eq:L:odd:sum:hyp} is true because both sides
  equal $h_1^{(k, k ; m)}$ by definition. If we subtract the identity
  \eqref{eq:L:odd:sum:hyp} from the corresponding identity with $d$ replaced
  by $d + 1$, we obtain
  \begin{eqnarray}
    &  & h_1^{(k + 2 d + 1, k + 2 d + 1 ; m - 4 d - 2)} + h_1^{(k + 2 d + 2,
    k + 2 d + 2 ; m - 4 d - 4)} \nonumber\\
    & = & h_1^{(k + 2 d + 2, k ; m - 2 d - 2)} - h_1^{(k + 2 d, k ; m - 2 d)}
    + \sum_{j = 0}^d h^{(k + 2 d, k + 2 j, m - 2 (d + j))} . 
    \label{eq:L:odd:induction}
  \end{eqnarray}
  By induction, the identity~\eqref{eq:L:odd:sum:hyp} is true for all $d
  \geq 0$ if we can show that \eqref{eq:L:odd:induction} holds for all $d
  \geq 0$.
  
  In the remainder, we prove \eqref{eq:L:odd:induction} by induction on $d$.
  The base case $d = 0$ of \eqref{eq:L:odd:induction} reads
  \begin{equation}
    h_1^{(k + 1, k + 1 ; m - 2)} + h_1^{(k + 2, k + 2 ; m - 4)} = h_1^{(k + 2,
    k ; m - 2)} - h_1^{(k, k ; m)} + h^{(k, k, m)} \label{eq:h3:h1kk}
  \end{equation}
  which can be proved in the same manner as Lemma~\ref{lem:3=2identity}.
  Alternatively, using that $h_1^{(k, k + 2 ; m - 2)} = - h_1^{(k + 1, k + 1 ;
  m - 2)}$, identity \eqref{eq:h3:h1kk} can be obtained as the special case
  $\ell = k$ of the identity \eqref{eq:h3:h1lk} which is proved below. Observe
  that replacing $d$ by $d - 1$ as well as $k$ by $k + 2$ and $m$ by $m - 4$
  in \eqref{eq:L:odd:induction} results in:
  \begin{eqnarray}
    &  & h_1^{(k + 2 d + 1, k + 2 d + 1 ; m - 4 d - 2)} + h_1^{(k + 2 d + 2,
    k + 2 d + 2 ; m - 4 d - 4)} \nonumber\\
    & = & h_1^{(k + 2 d + 2, k + 2 ; m - 2 d - 4)} - h_1^{(k + 2 d, k + 2 ; m
    - 2 d - 2)} + \sum_{j = 0}^{d - 1} h^{(k + 2 d, k + 2 j + 2, m - 2 (d + j
    + 1))}  \label{eq:L:odd:induction:2}
  \end{eqnarray}
  Let us fix a value $d > 0$ with the aim of proving
  \eqref{eq:L:odd:induction}. Since \eqref{eq:L:odd:induction:2} was obtained
  from \eqref{eq:L:odd:induction} by replacing $d$ with $d - 1$, we may assume
  as the induction hypothesis that \eqref{eq:L:odd:induction:2} holds for our
  fixed value of $d$. Note that the left-hand sides of
  \eqref{eq:L:odd:induction} and \eqref{eq:L:odd:induction:2} agree. Further
  note that the terms of the sum over $j$ in \eqref{eq:L:odd:induction:2} also
  appear in the corresponding sum in \eqref{eq:L:odd:induction}, except that
  the index $j$ is shifted and that \eqref{eq:L:odd:induction} has an
  additional term. Accordingly, subtracting \eqref{eq:L:odd:induction:2} from
  \eqref{eq:L:odd:induction} results in
  \begin{eqnarray*}
    0 & = & h_1^{(k + 2 d + 2, k ; m - 2 d - 2)} - h_1^{(k + 2 d, k ; m - 2
    d)} + h^{(k + 2 d, k, m - 2 d)}\\
    &  & - h_1^{(k + 2 d + 2, k + 2 ; m - 2 d - 4)} + h_1^{(k + 2 d, k + 2 ;
    m - 2 d - 2)} .
  \end{eqnarray*}
  Proving this fixed-length identity is equivalent to proving
  \eqref{eq:L:odd:induction}. If we replace $m$ by $m + 2 d$ and then replace
  $2 d$ by $\ell - k$ the identity takes the form
  \begin{equation}
    0 = h_1^{(\ell + 2, k ; m - 2)} - h_1^{(\ell, k ; m)} + h^{(\ell, k, m)} -
    h_1^{(\ell + 2, k + 2 ; m - 4)} + h_1^{(\ell, k + 2 ; m - 2)} .
    \label{eq:h3:h1lk}
  \end{equation}
  After dividing by $h^{(\ell, k, m)}$, \eqref{eq:h3:h1lk} is equivalent to
  \begin{equation}
    \frac{h_1^{(\ell, k ; m)}}{h^{(\ell, k, m)}} - \frac{h_1^{(\ell + 2, k ; m
    - 2)}}{h^{(\ell, k, m)}} - \frac{h_1^{(\ell, k + 2 ; m - 2)}}{h^{(\ell, k,
    m)}} + \frac{h_1^{(\ell + 2, k + 2 ; m - 4)}}{h^{(\ell, k, m)}} = 1.
    \label{eq:h3:h1lk:rat}
  \end{equation}
  By the formulas defining $h$ and $h_1$, each of the four terms on the
  left-hand side is a rational function in $\ell, k$ and $m$. For instance,
  the first term is
  \begin{equation*}
    \frac{h_1^{(\ell, k ; m)}}{h^{(\ell, k, m)}} = \frac{(\ell + 1) (\ell +
     2) k (k + 1)}{(k + m) (\ell + m + 1) (\ell - m + 2) (k - m + 1)} .
  \end{equation*}
  Using these rational functions in \eqref{eq:h3:h1lk:rat}, the identity
  \eqref{eq:h3:h1lk:rat} is readily confirmed.
\end{proof}

Since it appears to have independent value, we record the following which is a
special case of the identity \eqref{eq:h3:h1lk} that was used to complete the
proof of Theorem~\ref{thm:L:odd:sum}. We further note that
Lemma~\ref{lem:3=2identity} can also be obtained as a special case of
\eqref{eq:h3:h1lk} by setting $\ell = k$ and observing that $h_1^{(k, k + 2 ;
m - 2)} = - h_1^{(k + 1, k + 1 ; m - 2)}$.

\begin{corollary}
  \label{cor:f3:f1lk}For integers $k \geq m \geq 4$ and $\ell
  \geq k + 2$,
  \begin{equation}
    f^{(\ell, k, m)} = f^{(\ell, k, 1^m)} - f^{(\ell, k + 2, 1^{m - 2})} -
    f^{(\ell + 2, k, 1^{m - 2})} + f^{(\ell + 2, k + 2, 1^{m - 4})} .
  \end{equation}
\end{corollary}

Theorem~\ref{thm:L:odd:sum} provides an evaluation of $L_{2 d + 1} (k, m)$
when $m \leq k$. The analytic perspective adopted in the proof of
Theorem~\ref{thm:L:odd:sum} allows us to obtain a corresponding evaluation
when $m$ is large.

\begin{corollary}
  \label{cor:L:odd:sum:2}Let $d \geq 0$ as well as $k, m \geq 2$ be
  integers. If $m \geq k + 6 d - 3$ then
  \begin{equation}
    L_{2 d + 1} (k, m) = f^{(k + 2 d, k, 1^{m - 2 d})} + \sum_{r = 0}^{d - 1}
    \sum_{j = 0}^r f^{(m - 2 (r + j + 1), k + 2 r + 1, k + 2 j + 1)} .
    \label{eq:L:odd:sum:2}
  \end{equation}
\end{corollary}

\begin{proof}
  From the definition of $h^{(x, y, z)}$ in the proof of
  Theorem~\ref{thm:L:odd:sum} it is readily observed that
  \begin{equation}
    h^{(x, y, z)} = h^{(z - 2, x + 1, y + 1)} . \label{eq:h:swap}
  \end{equation}
  In particular,
  \begin{equation*}
    h^{(k + 2 r, k + 2 j, m - 2 (r + j))} = h^{(m - 2 (r + j) - 2, k + 2 r +
     1, k + 2 j + 1)} .
  \end{equation*}
  Applying this relation to the identity \eqref{eq:L:odd:sum:hyp} results in
  \begin{equation}
    \sum_{j = 0}^{2 d} h_1^{(k + j, k + j ; m - 2 j)} = h_1^{(k + 2 d, k ; m -
    2 d)} + \sum_{r = 0}^{d - 1} \sum_{j = 0}^r h^{(m - 2 (r + j + 1), k + 2 r
    + 1, k + 2 j + 1)} . \label{eq:L:odd:sum:2:h}
  \end{equation}
  The desired identity \eqref{eq:L:odd:sum:2} follows immediately from
  \eqref{eq:L:odd:sum:2:h} provided that each term in \eqref{eq:L:odd:sum:2:h}
  involving $h$ is of the form $h^{(x, y, z)}$ where $(x, y, z)$ is a
  partition (and likewise for $h_1$). This is ensured by the condition $m
  \geq k + 6 d - 3$ which implies that the part $m - 2 (r + j + 1)$ is no
  less than the part $k + 2 r + 1$ for all relevant choices of $r$ and $j$.
\end{proof}

In the case $d = 1$, Corollary~\ref{cor:L:odd:sum:2} reduces to the second
case of Lemma~\ref{lem:3=2identity} while the case $d = 2$ results in the
identity \eqref{eq:L5:2}. The intermediate identities discussed in
Example~\ref{eg:L5} can be obtained by inspecting each term of the
corresponding analytic identity \eqref{eq:L:odd:sum:hyp}.

As an application, we use the results of this section to give an alternative
proof of the part of Theorem~\ref{thm: Generalizedknapsack} that refines
Theorem~\ref{thm:equalitySpecht}.

\begin{theorem}
  Let $n$ and $k$ be positive integers such that $k \le \tfrac{n}{2} - 1$.
  If $k \leq \lceil n / 3 \rceil$ or if $n \equiv k \pmod{2}$, then
  \begin{equation}
    \sum_{\lambda \in X_0 (n, k)} f^{\lambda} = f^{(k, k, 1^{n - 2 k})} +
    f^{(k + 1, k + 1, 1^{n - 2 k - 2})} \label{eq:L2:X1}
  \end{equation}
  where $X_0 (n, k) = \left\{ (\lambda_1, k, \lambda_3) \vdash n \, | \, k
  \equiv \lambda_3 \pmod{2} \right\}$ as in
  Theorem~\ref{thm: Generalizedknapsack}.
\end{theorem}

\begin{proof}
  Set $m = n - 2 k$. Upon subtracting \eqref{eq:L:odd:sum:2:h} from
  \eqref{eq:L:odd:sum:2:h} with $d$ replaced by $d + 1$, we obtain
  \begin{eqnarray*}
    &  & h_1^{(k + 2 d + 1, k + 2 d + 1 ; m - 4 d - 2)} + h_1^{(k + 2 d + 2,
    k + 2 d + 2 ; m - 4 d - 4)}\\
    & = & h_1^{(k + 2 d + 2, k ; m - 2 d - 2)} - h_1^{(k + 2 d, k ; m - 2 d)}
    + \sum_{j = 0}^d h^{(m - 2 (d + j + 1), k + 2 d + 1, k + 2 j + 1)} .
  \end{eqnarray*}
  Replacing $k$ by $k - 2 d - 1$ and $m$ by $m + 4 d + 2$, this takes the form
  \begin{eqnarray}
    &  & h_1^{(k, k ; m)} + h_1^{(k + 1, k + 1 ; m - 2)}  \label{eq:L2:h:d}\\
    & = & h_1^{(k + 1, k - 2 d - 1 ; m + 2 d)} - h_1^{(k - 1, k - 2 d - 1 ; m
    + 2 d + 2)} + \sum_{j = 0}^d h^{(m + 2 j, k, k - 2 j)}, \nonumber
  \end{eqnarray}
  where we replaced $j$ by $d - j$ in the final sum. By definition,
  \begin{equation*}
    h_1^{(k + 1, k - 2 d - 1 ; m + 2 d)} = \frac{(m + 2 k) !}{(k + 1) ! (k -
     2 d - 2) ! (m + 2 d) !} \frac{2 d + 3}{(m + k + 2 d + 2) (m + k - 1)}
  \end{equation*}
  and, because of the factor $(k - 2 d - 2) !$ in the denominator, this
  expression vanishes for either $d = k / 2$ or $d = (k - 1) / 2$, depending on whether $k$ is even or odd. Likewise,
  the term $h_1^{(k - 1, k - 2 d - 1 ; m + 2 d + 2)}$ vanishes for $d =
  \lfloor k / 2 \rfloor$. Hence, specializing \eqref{eq:L2:h:d} to $d =
  \lfloor k / 2 \rfloor$ results in
  \begin{equation}
    h_1^{(k, k ; m)} + h_1^{(k + 1, k + 1 ; m - 2)} = \sum_{j = 0}^{\lfloor k
    / 2 \rfloor} h^{(m + 2 j, k, k - 2 j)} . \label{eq:L2:h}
  \end{equation}
  If $m \geq k$ or, equivalently, $3 k \leq n$ then the summands on
  the right-hand side are $f^{(m + 2 j, k, k - 2 j)}$ and the sum equals the
  sum over $f^{\lambda}$ for $\lambda \in X_0 (n, k)$. This proves
  \eqref{eq:L2:X1} in the case $k \leq n / 3$. To show that
  \eqref{eq:L2:X1} continues to hold for $k \leq \lceil n / 3 \rceil$ or,
  equivalently, $3 k \leq n + 2$, we need to also consider the cases $m =
  k - 1$ and $m = k - 2$. In those cases, the right-hand side of
  \eqref{eq:L2:h} equals $h^{(m, k, k)}$ plus the sum over $f^{\lambda}$ for
  $\lambda \in X_0 (n, k)$. In other words, if $m = k - 1$ or $m = k - 2$ then
  \begin{equation*}
    f^{(k, k, 1^m)} + f^{(k + 1, k + 1, 1^{m - 2})} = h^{(m, k, k)} +
     \sum_{\lambda \in X_0 (n, k)} f^{\lambda} .
  \end{equation*}
  By the definition of $h$, we see that
  \begin{equation*}
    h^{(m, k, k)} = \frac{(m + 2 k) !}{(m + 2) ! (k + 1) !k!} \; (m - k + 1)
     (m - k + 2)
  \end{equation*}
  vanishes for $m = k - 1$ as well as $m = k - 2$. This completes the proof of
  \eqref{eq:L2:X1} in the case $k \leq \lceil n / 3 \rceil$.
  
  In the remainder, we therefore assume that $3 k > n$ and $n \equiv k
  \pmod{2}$. Equivalently, in terms of $m$, we assume $k >
  m$ and $m \equiv k \pmod{2}$. Note that the summand $s
  (j) = h^{(m + 2 j, k, k - 2 j)}$ of \eqref{eq:L2:h} satisfies the symmetry
  \begin{equation*}
    s \left(\frac{k - m}{2} - 1 - j \right) = - s (j)
  \end{equation*}
  which follows immediately from the representation
  \begin{equation*}
    s (j) = \frac{(m + 2 k) ! (m + 2 j - k + 1) (m - k + 4 j + 2) (2 j +
     1)}{(m + 2 j + 2) ! (k + 1) ! (k - 2 j) !} .
  \end{equation*}
  Consequently,
  \begin{equation*}
    \sum_{j = 0}^{(k - m) / 2 - 1} s (j) = \sum_{j = 0}^{(k - m) / 2 - 1} s
     \left(\frac{k - m}{2} - 1 - j \right) = - \sum_{j = 0}^{(k - m) / 2 - 1}
     s (j)
  \end{equation*}
  so that these sums all vanish. Combined with \eqref{eq:L2:h}, we therefore
  have
  \begin{equation*}
    f^{(k, k, 1^m)} + f^{(k + 1, k + 1, 1^{m - 2})} = \sum_{j = (k - m) /
     2}^{\lfloor k / 2 \rfloor} h^{(m + 2 j, k, k - 2 j)} = \sum_{\lambda \in
     X_0 (n, k)} f^{\lambda}
  \end{equation*}
  as claimed.
\end{proof}

The following examples illustrate the two different cases where
\eqref{eq:L2:h} has additional terms compared to \eqref{eq:L2:X1}.

\begin{example}
  Let $n = 13$ and $k = 5$ so that $k > n / 3$ but $k \leq \lceil n / 3
  \rceil$. Then \eqref{eq:L2:h} takes the form
  \begin{equation*}
    f^{(5, 5, 1^3)} + f^{(6, 6, 1^1)} = h^{(3, 5, 5)} + h^{(5, 5, 3)} +
     h^{(7, 5, 1)} .
  \end{equation*}
  Since $h^{(3, 5, 5)} = 0$, this reduces to $f^{(5, 5, 1^3)} + f^{(6, 6,
  1^1)} = f^{(5, 5, 3)} + f^{(7, 5, 1)}$ which is \eqref{eq:L2:X1}.
\end{example}

\begin{example}
  To illustrate the case $k > n / 3$ and $n \equiv k \pmod{2}$, let $n = 20$ and $k = 8$. In that case the identity \eqref{eq:L2:h}
  becomes
  \begin{eqnarray*}
    f^{(8, 8, 1^4)} + f^{(9, 9, 1^2)} & = & \sum_{j = 0}^4 h^{(4 + 2 j, 8, 8 -
    2 j)}\\
    & = & h^{(4, 8, 8)} + h^{(6, 8, 6)} + f^{(8, 8, 4)} + f^{(10, 8, 2)} +
    f^{(12, 8)} .
  \end{eqnarray*}
  The terms $h^{(4, 8, 8)} = 1385670$ and $h^{(6, 8, 6)} = - 1385670$ cancel
  and we arrive at
  \begin{equation*}
    f^{(8, 8, 1^4)} + f^{(9, 9, 1^2)} = f^{(8, 8, 4)} + f^{(10, 8, 2)} +
     f^{(12, 8)}
  \end{equation*}
  as predicted by \eqref{eq:L2:X1}.
\end{example}

It is natural to wonder if the results in this section extend to the sums $L_d
(k, m)$ when $d$ is even. For general values of $k$ and $m$, our numerical
data suggests that no analogous identities exist. However, more specialized
formulas appear to hold when $k$ and $m$ have opposite parity. We do not
pursue this question here.

\section{Concluding remarks}

Our main results, Theorem~\ref{thm: Generalizedknapsack} and
Theorem~\ref{thm:L:odd:sum}, as well as auxiliary results such as Lemma~\ref{lem:3=2identity}, are
all instances of identities between character degrees of the form
\begin{equation}
  \sum_{\lambda \in X_n} f^{\lambda} = \sum_{\mu \in Y_n} f^{\mu}
  \label{eq:sum:XY}
\end{equation}
where $X_n$ and $Y_n$ are certain sets of partitions of size $n$. One can
think of \eqref{eq:sum:XY} as packing two ``knapsacks'' with the various
$\lambda$ as books and the $f^{\lambda}$ as the corresponding weights so that
both knapsacks are filled to equal weight. It is natural to wonder if one can
find further interesting results of the form \eqref{eq:sum:XY} for natural
sets $X_n$ and $Y_n$ of partitions.

There is only one place we have seen results of the form \eqref{eq:sum:XY}, namely in the context of modular
representation theory. Starting with a partition $\mu$ and an integer $k$, if
we add a single rim $k$-hook to $\mu$ in every possible way, we get an
alternating sum of characters that vanishes on any permutation not containing
a $k$-cycle. In particular, it vanishes on the identity permutation, so that
we obtain an identity on the character degrees. For details, we refer to
\cite[21.7]{James_1984}. For an example, if $\mu = (3, 1)$ and $k = 6$ then
we obtain
\begin{equation*}
  f^{(9, 1)} - f^{(6, 4)} + f^{(4, 4, 2)} - f^{(3, 3, 2, 1^2)} + f^{(3, 2, 2,
   1^3)} - f^{(3, 1^7)} = 0.
\end{equation*}
We pose the very speculative question: can our identities be related to those arising from this hook-wrapping
procedure? 

We remark that a variation of \eqref{eq:sum:XY} appears in the work of Regev
\cite{regev-syt-kl-hooks} on the number of standard Young tableaux contained
in the $(k, \ell)$-hook, which are tableaux whose shape is contained in
\begin{equation*}
  H (k, \ell ; n) = \left\{ \lambda = (\lambda_1, \lambda_2, \ldots) :
   \lambda \vdash n \text{ and } \lambda_{k + 1} \leq \ell \right\} .
\end{equation*}
Specifically, it is shown in \cite[Proposition~3.2]{regev-syt-kl-hooks} that
\begin{equation}
  \sum_{\lambda \in H^{\ast} (2, 2 ; 2 m)} f^{\lambda} = \sum_{\lambda \in H
  (4, 0 ; 2 m - 2)} f^{\lambda} \label{eq:sum:XY:regev}
\end{equation}
where the partitions on the left-hand side are
\begin{equation*}
  H^{\ast} (2, 2 ; 2 m) = \{ (k, k, 2^{m - k}) \vdash 2 m : k \in \{ 2,
   \ldots, m \} \}
\end{equation*}
while the right-hand sum is over partitions of size $2 m - 2$ with up to four
parts. This differs from \eqref{eq:sum:XY} in that the sums in
\eqref{eq:sum:XY:regev} are over partitions of differing size. We note that
both sides of \eqref{eq:sum:XY:regev} are counted by the product $C_{m - 1}
C_m$ of Catalan numbers. This follows from the result of Gouyou-Beauchamps
\cite{gouyou-beauchamps-syt-height45} that the number of standard Young
tableaux of size $n$ with at most four rows is $C_{\lfloor (n + 1) / 2
\rfloor} C_{\lceil (n + 1) / 2 \rceil}$.

Finally, we wonder if our results, especially Theorem~\ref{thm:
Generalizedknapsack} and Theorem~\ref{thm:L:odd:sum}, have a natural
representation-theoretic interpretation. In a similar direction, can one find
combinatorial proofs of our identities via bijections on the corresponding
sets of standard Young tableaux? 
Can one find a bijective proof of Theorem \ref{thm: Generalizedknapsack}? 
Echoing a problem already raised at the end of \cite{HemmerStraubWestrem2025}, even for the unrefined version
of Theorem~\ref{thm: Generalizedknapsack}, namely Theorem~\ref{thm:equalitySpecht}, we do not know of
a bijective proof (but see footnote \footref{fn}). One would need a bijection between, on one side, all
standard Young tableaux of shape $\lambda$ where $\lambda$ is a partition of
$n$ into three parts of equal parity and, on the other side, standard Young
tableaux of shape $(k, k, 1^{n - 2 k})$ where $k$ varies. A bijection between
the latter and Riordan paths is known (see, for instance,
\cite{HemmerStraubWestrem2025}), and bijections are known between Young
tableaux with up to three rows and Motzkin paths (see, for instance,
\cite{gmtw-dyck-syt} or \cite{MatsakisMotzkinInspired}). It is therefore
natural to look for a bijection between Motzkin paths and Young tableaux with
up to three rows that sends Riordan paths to Young tableaux of shape $\lambda$
where $\lambda$ has parts of equal parity.

{\textbf{Acknowledgements. }} We would like to thank the anonymous referee for a very detailed and thorough reading of an earlier draft of this paper.

\appendix

\section{A complete branching proof of Theorem \ref{thm: Generalizedknapsack}}
\label{sec:complete-branching-proof}

In this appendix we give the details of the simultaneous induction proof of
Theorem~\ref{thm: Generalizedknapsack}.  The notation is designed to keep the case analysis
short while still making clear that all cases are covered. The proof will be organized based on the size of $k$ and whether or not it is congruent to $n$ modulo 2.

To simplify the notation and to match Lemma \ref{lem:3=2identity}, throughout this appendix, put
\[
        m:=n-2k.
\]
Thus $n\equiv m\pmod 2$, and the condition
$k>\lceil n/3\rceil$ is equivalent to $k\ge m+3$.  

We use the convention that any displayed term $f^\lambda$ is zero if
$\lambda$ is not a partition; for example, a term involving $1^t$ with
$t<0$ is omitted. That way, the statement of Theorem~\ref{thm: Generalizedknapsack} holds true for all $n\geq 3$ and $k \geq 0$. The theorem is readily checked for all values of $n \le 10$.
For example, for $n=10$, the identities \eqref{eq:firstdimension} and
\eqref{eq:seconddimension} for $k=1,2,3,4$ read:
\begin{align*}
f^{(8,1,1)} &= f^{(1,1,1^{8})} + f^{(2,2,1^{6})} = 36, \\
f^{(9,1)} &= f^{(2,1,1^{7})} = 9, \\[4pt]
f^{(8,2)} + f^{(6,2,2)} &= f^{(2,2,1^{6})} + f^{(3,3,1^{4})} = 260, \\
f^{(7,2,1)} &= f^{(3,2,1^{5})} = 160, \\[4pt]
f^{(6,3,1)} + f^{(4,3,3)} &= f^{(3,3,1^{4})} + f^{(4,4,1^{2})} = 525, \\
f^{(7,3)} + f^{(5,3,2)} &= f^{(4,3,1^{3})} = 525, \\[4pt]
f^{(6,4)} + f^{(4,4,2)} &= f^{(4,4,1^{2})} + f^{(5,5)} = 342, \\
f^{(5,4,1)} &= f^{(5,4,1)} = 288.
\end{align*}
When $k$ or $m$ is very small, specifically, when
$k\le1$ or $m\le2$, the inductive argument
below would use Lemma~\ref{lem:3=2identity} with parameters outside its
hypotheses. In these cases (which occur for arbitrarily large $n$) we will verify the main theorem directly.

Next, we consolidate the two sums in Theorem \ref{thm: Generalizedknapsack}. For $e\in\{0,1\}$, define
\[
S_e(n,k)
 = \sum_{\lambda \in X_e(n,k)} f^\lambda
 =
 \sum_{\substack{0\le r\le \min(k,m)\\ r\equiv k+e\pmod 2}}
 f^{(n-k-r,k,r)}.
\]
Further, define
\[
A(n,k)=f^{(k,k,1^m)}+f^{(k+1,k+1,1^{m-2})}
\]
and
\[
B(n,k)=f^{(k+1,k,1^{m-1})}.
\]
With this notation, Theorem~\ref{thm: Generalizedknapsack} says that
\[
S_0(n,k)=A(n,k),\qquad S_1(n,k)=B(n,k)
\]
unless $k\ge m+3$ and $k\not\equiv m\pmod 2$, in which case the two roles are
switched:
\[
S_0(n,k)=B(n,k),\qquad S_1(n,k)=A(n,k).
\]

Equivalently, define
\[
\sigma(n,k)=
\begin{cases}
1, & \text{if } k\ge m+3 \text{ and } k\not\equiv m\pmod 2,\\
0, & \text{otherwise.}
\end{cases}
\]
Then the theorem asserts
\begin{equation}
\label{eq:appendix-main-claim}
        S_{\sigma(n,k)}(n,k)=A(n,k),
        \qquad
        S_{1-\sigma(n,k)}(n,k)=B(n,k).
\end{equation}
Before considering the general case, we first verify
\eqref{eq:appendix-main-claim} when  $k\le1$ or $m\le2$.

First suppose $k\le1$, so that $\sigma(n,k)=0$. For $k=0$, we have
$X_0(n,0)=\{(n)\}$ and $X_1(n,0)=\varnothing$, and the claims reduce to
$f^{(n)}=f^{(1^n)}$, both equal to $1$, and to $0=0$. For
$k=1$, we have $X_0(n,1)=\{(n-2,1,1)\}$ and $X_1(n,1)=\{(n-1,1)\}$, and the
claims become
\[
  f^{(n-2,1,1)}=f^{(1,1,1^{n-2})}+f^{(2,2,1^{n-4})},
  \qquad
  f^{(n-1,1)}=f^{(2,1,1^{n-3})}.
\]
By Lemma~\ref{lem:closedformab111}, both sides of the first identity equal
$\binom{n-1}{2}$ and both sides of the second equal $n-1$.

Now suppose $m\le2$. If $k<m+3$ then $n = m+2k \le 10$, and we already verified these base cases. We may therefore assume that $k\ge m+3$. In that case, $\sigma(n,k)=1$ exactly when $m$ is odd. 

For $m=0$, both claims reduce to
$f^{(k,k)}=f^{(k,k)}$, all other displayed terms being omitted. For $m=1$,
the sets are $X_0(n,k)=\{(k+1,k)\}$ and $X_1(n,k)=\{(k,k,1)\}$, and the
claims reduce to $f^{(k,k,1)}=f^{(k,k,1)}$ and
$f^{(k+1,k)}=f^{(k+1,k)}$. For $m=2$, we have
$X_0(n,k)=\{(k+2,k),(k,k,2)\}$ and $X_1(n,k)=\{(k+1,k,1)\}$; the second
claim is  $f^{(k+1,k,1)}=f^{(k+1,k,1)}$, while the first reads
\[
  f^{(k+2,k)}+f^{(k,k,2)}=f^{(k,k,1^2)}+f^{(k+1,k+1)}.
\]
This follows from Lemma~\ref{lem:closedformab111}: both sides equal
\[
\frac{k^2+3}{k+3}\binom{2k+1}{k}.
\]

In the induction below we may therefore assume $k\ge2$ and $m\ge3$. In this
range, Lemma~\ref{lem:3=2identity} and
Proposition~\ref{prop: boundaryspecialcase} are only invoked with
parameters within their hypotheses, and every displayed term is a genuine
character degree; the convention above is then needed only for checking the
finitely many base cases of small $n$.

So we are left to prove the two equations in \eqref{eq:appendix-main-claim} by simultaneous induction on $n$.
The base cases for small $n$ are immediate.  Assume then that
\eqref{eq:appendix-main-claim} holds for all smaller values of $n$.

Let $\br$ denote the result of applying the branching rule to a sum of
character degrees.  In particular, for a partition $\lambda\vdash n$,
\[
\br(f^\lambda)
=
\sum_{\mu} f^\mu,
\]
where $\mu$ runs over the partitions obtained from $\lambda$ by removing one
removable node.  The branching rule says that
\[
        f^\lambda=\br(f^\lambda).
\]
Consequently, it is enough to compare the branched forms of the two sides of
each desired identity, which will let us apply the inductive hypotheses.

For convenience, set
\[
\begin{aligned}
A_0&=A(n-1,k)
   =f^{(k,k,1^{m-1})}+f^{(k+1,k+1,1^{m-3})},\\
B_0&=B(n-1,k)
   =f^{(k+1,k,1^{m-2})},\\
A_1&=A(n-1,k-1)
   =f^{(k-1,k-1,1^{m+1})}+f^{(k,k,1^{m-1})},\\
B_1&=B(n-1,k-1)
   =f^{(k,k-1,1^m)}.
\end{aligned}
\]
Applying the branching rule directly to $A(n,k)$ gives four terms, which one checks to be
\begin{equation}
\label{eq:branch-A}
        \br(A(n,k))=A_0+B_1+B_0.
\end{equation}
Likewise,
\begin{equation}
\label{eq:branch-B}
\br(B(n,k))
 =
 f^{(k,k,1^{m-1})}
 +f^{(k+1,k-1,1^{m-1})}
 +f^{(k+1,k,1^{m-2})}.
\end{equation}

We shall use Lemma~\ref{lem:3=2identity} in the following rearranged form,
obtained by replacing its parameters $k,m$ with $k-1,m+1$:
\begin{equation}
\label{eq:appendix-rearranged-L}
\begin{aligned}
f^{(k+1,k-1,1^{m-1})}
&=
f^{(k-1,k-1,1^{m+1})}
+f^{(k,k,1^{m-1})}
+f^{(k+1,k+1,1^{m-3})}
-E,
\end{aligned}
\end{equation}
where
\[
E=
\begin{cases}
f^{(k-1,k-1,m+1)}, & m\le k-2,\\
0, & m=k-1 \text{ or } m=k,\\
f^{(m-1,k,k)}, & m\ge k+1.
\end{cases}
\]
Combining \eqref{eq:appendix-rearranged-L} with the branching formula \eqref{eq:branch-B} for
$B(n,k)$, we get
\begin{equation}
\label{eq:branch-BwithE}
        \br(B(n,k))=B_0+A_0+A_1-E.
\end{equation}

So in \eqref{eq:branch-A} and \eqref{eq:branch-BwithE} we have the branched versions of the right-hand sides of \eqref{eq:appendix-main-claim}. We next determine the branching behavior of the left-hand sides. To compute $\br(S_e(n,k))$ (see \eqref{eq:branch-Se} for the resulting formula), we remove boxes row by row in the respective Young diagrams. Keep in mind that a box is removable from the first (resp.\ second) row of $\lambda$ only if $\lambda_1>\lambda_2$ (resp.\ $\lambda_2>\lambda_3$), and from the third row only if $\lambda_3\ge1$.

Consider a diagram in $X_e(n,k)$ (here, and in what follows, we identify partitions with their Young diagrams). Removing a box from the first row of the diagram does not affect the parity condition on
the second and third parts, and so produces a diagram in $X_e(n-1,k)$. (Note that the only diagram in $X_e(n,k)$ with no removable
first-row box is $(k,k,m)$, which occurs when $m\le k$ and
$m\equiv k+e\pmod 2$.) Conversely, each diagram in
$X_e(n-1,k)$ arises from exactly one diagram in $X_e(n,k)$ by adding a box to the first row. As such, we obtain exactly $S_e(n-1,k)$ as the contribution to $\br(S_e(n,k))$ arising from removing a box from the first row.

Removing a box from the second row changes $k$ to $k-1$ and reverses the
parity condition, producing diagrams in $X_{1-e}(n-1,k-1)$. (Note that the only diagram in $X_e(n,k)$ with no removable second-row box is $(m,k,k)$, which occurs when $m\ge k$ and $e=0$.) Conversely, with at most one exception, each diagram in $X_{1-e}(n-1,k-1)$ arises from exactly one diagram in $X_e(n,k)$ by adding a box to the second row. The potential exception is
\[
        (k-1,k-1,m+1),
\]
which would give rise to the invalid diagram $(k-1,k,m+1)$ upon adding to the second row. This diagram is in $X_{1-e}(n-1,k-1)$, and hence missing, precisely when $m\le k-2$ and $m+1\equiv k+e\pmod 2$.

Lastly, removing a box from the third row again reverses the parity
condition and produces diagrams in $X_{1-e}(n-1,k)$. (Note that the only diagram in $X_e(n,k)$ with no removable third-row box is $(n-k,k)$, which occurs when $k\equiv k+e \pmod 2$.) Conversely, with again at most one exception, each diagram in $X_{1-e}(n-1,k)$ arises from exactly one diagram in $X_e(n,k)$ by adding a box to the third row. The potential exception is
\[
        (m-1,k,k),
\]
which would give rise to the invalid diagram $(m-1,k,k+1)$ upon adding to the third row. This diagram is in $X_{1-e}(n-1,k)$, and hence missing, precisely when $e=1$ and $m\ge k+1$.

Combining these three cases, we therefore find
\begin{equation}
\label{eq:branch-Se}
\begin{aligned}
\br(S_e(n,k))
&=
S_e(n-1,k)
+S_{1-e}(n-1,k-1)
+S_{1-e}(n-1,k)\\
&\quad
-\alpha_e f^{(k-1,k-1,m+1)}
-\beta_e f^{(m-1,k,k)},
\end{aligned}
\end{equation}
where
\[
\alpha_e=
\begin{cases}
1, & m\le k-2 \text{ and } m+1\equiv k+e\pmod 2,\\
0, & \text{otherwise,}
\end{cases}
\]
and
\[
\beta_e=
\begin{cases}
1, & e=1 \text{ and } m\ge k+1,\\
0, & \text{otherwise.}
\end{cases}
\]

To complete the proof, we will apply the inductive hypotheses to the first three terms on the right-hand side of \eqref{eq:branch-Se} and confirm, together with extra terms whenever $\alpha_e$ and/or $\beta_e$ is one, that the result matches \eqref{eq:appendix-main-claim} branched using \eqref{eq:branch-A} and \eqref{eq:branch-BwithE}. 

We divide the proof into four cases.

\begin{description}[leftmargin=0pt, labelindent=0pt, itemsep=\smallskipamount]

\item[Case 1: $k\le m+1$]

Since $k\le m+1<m+3$, we have
\[
        \sigma(n,k)=0.
\]
Moreover, both lower pairs $(n-1,k)$ and $(n-1,k-1)$ are also in the
unswitched regime, meaning that
\[
        \sigma(n-1,k)=0,\qquad \sigma(n-1,k-1)=0.
\]

First consider $S_0(n,k)$, that is $e=0$.  In that case, we have $\alpha_e=\beta_e=0$ in \eqref{eq:branch-Se}.  Hence, by the induction hypothesis applied to the right-hand side of \eqref{eq:branch-Se}, we find
\[
\br(S_0(n,k))
 =
A_0+B_1+B_0.
\]
By \eqref{eq:branch-A}, this is exactly $\br(A(n,k))$.  This proves that
$S_0(n,k)=A(n,k)$.

Now consider $S_1(n,k)$, that is $e=1$.  In that case, $\alpha_e=0$. On the other hand, $\beta_e=1$ precisely when $m\ge k+1$, in which case the term $-f^{(m-1,k,k)}$ is present on the right-hand side of \eqref{eq:branch-Se}.  Thus, using the induction hypothesis,
\[
\br(S_1(n,k))
 =
B_0+A_1+A_0-E,
\]
where $E=0$ if $m=k-1$ or $m=k$, and
$E=f^{(m-1,k,k)}$ if $m\ge k+1$.  By \eqref{eq:branch-BwithE}, this is
$\br(B(n,k))$.  Hence, $S_1(n,k)=B(n,k)$. This completes the proof in Case 1.

\item[Case 2: $k\ge m+2$ and $k\equiv m\pmod 2$]

Here $k$ and $n$ have the same parity, so we are still in the
unswitched regime:
\[
        \sigma(n,k)=0.
\]
For the lower pairs, one checks that
\[
        \sigma(n-1,k)=1,
        \qquad
        \sigma(n-1,k-1)=0.
\]

For $S_0(n,k)$, we have $\alpha_e=\beta_e=0$ in \eqref{eq:branch-Se}.  Therefore, the induction hypothesis gives
\[
\br(S_0(n,k))
 =
B_0+B_1+A_0.
\]
By \eqref{eq:branch-A}, this equals $\br(A(n,k))$.  Hence,
$S_0(n,k)=A(n,k)$.

For $S_1(n,k)$, we have $\alpha_e=1$ and $\beta_e=0$. Accordingly,
\[
\br(S_1(n,k))
 =
A_0+A_1+B_0-f^{(k-1,k-1,m+1)}.
\]
Since $m\le k-2$, the final term exactly matches the
term $E$ in \eqref{eq:branch-BwithE}. Hence, $\br(S_1(n,k)) = \br(B(n,k))$, showing that $S_1(n,k)=B(n,k)$. This proves the theorem in Case 2.

\item[Case 3: $k=m+3$]

In this case $k\not\equiv m\pmod 2$ and $k\ge m+3$, so we are in the switched regime:
\[
        \sigma(n,k)=1.
\]
On the other hand, the two lower pairs are both unswitched:
\[
        \sigma(n-1,k)=0,
        \qquad
        \sigma(n-1,k-1)=0.
\]

Since $k-1=(m+1)+1$, it follows from Proposition~\ref{prop: boundaryspecialcase} that
\[
f^{(k-1,k-1,1^{m+1})}+f^{(k,k,1^{m-1})}
=
f^{(k,k-1,1^m)}.
\]
Equivalently, this means that, in the present case,
\begin{equation}
\label{eq:appendix-boundary-use}
        A_1=B_1.
\end{equation}

Starting with $S_1(n,k)$, we have $\alpha_e=\beta_e=0$ in \eqref{eq:branch-Se}. Therefore, using the induction hypothesis and then \eqref{eq:appendix-boundary-use},
\[
\br(S_1(n,k))
 =
B_0+A_1+A_0
 =
B_0+B_1+A_0.
\]
By \eqref{eq:branch-A}, this equals $\br(A(n,k))$.  Hence,
$S_1(n,k)=A(n,k)$.

Now consider $S_0(n,k)$, in which case we have $\alpha_e=1$ and $\beta_e=0$ in \eqref{eq:branch-Se}. Again, using the induction hypothesis and then \eqref{eq:appendix-boundary-use}, we find
\[
\br(S_0(n,k))
 =
A_0+B_1+B_0-f^{(k-1,k-1,m+1)}
 =
A_0+A_1+B_0-f^{(k-1,k-1,m+1)}.
\]
Since $m\le k-2$, the final term exactly matches the term $E$ in \eqref{eq:branch-BwithE}. Hence, $\br(S_0(n,k)) = \br(B(n,k))$, showing that $S_0(n,k)=B(n,k)$. This completes the proof of the theorem in Case 3.

\item[Case 4: $k\ge m+5$ and $k\not\equiv m\pmod 2$]

This is the remaining switched case:
\[
        \sigma(n,k)=1.
\]
For the lower pairs we have
\[
        \sigma(n-1,k)=0,
        \qquad
        \sigma(n-1,k-1)=1.
\]

Starting with $S_1(n,k)$, we have $\alpha_e=\beta_e=0$ in \eqref{eq:branch-Se}. Therefore, the induction hypothesis gives
\[
\br(S_1(n,k))
 =
B_0+B_1+A_0.
\]
By \eqref{eq:branch-A}, this is $\br(A(n,k))$. Hence, $S_1(n,k)=A(n,k)$.

Finally, consider $S_0(n,k)$, in which case we have $\alpha_e=1$ and $\beta_e=0$ in \eqref{eq:branch-Se}. Using the induction hypothesis, we obtain
\[
\br(S_0(n,k))
 =
A_0+A_1+B_0-f^{(k-1,k-1,m+1)}.
\]
Since $m\le k-5$, the final term matches the term $E$ in \eqref{eq:branch-BwithE}. This shows that $\br(S_0(n,k)) = \br(B(n,k))$. Hence, $S_0(n,k)=B(n,k)$. This proves the theorem in Case 4.
\end{description}

The four cases above exhaust all possibilities for $k$ relative to $m$. The simultaneous induction is therefore complete, and Theorem~\ref{thm: Generalizedknapsack} follows.

\bibliographystyle{alpha}
\bibliography{hookreferences}

\end{document}